\pdfoutput=1
\RequirePackage{ifpdf}
\ifpdf % We are running pdfTeX in pdf mode
\documentclass[pdftex]{sigma}
\else
\documentclass{sigma}
\fi

\begin{document}

\allowdisplaybreaks

\renewcommand{\PaperNumber}{050}

\FirstPageHeading

\ShortArticleName{A Connection Formula for the $q$-Conf\/luent Hypergeometric Function}

\ArticleName{A Connection Formula
\\
for the $\boldsymbol{q}$-Conf\/luent Hypergeometric Function}

\Author{Takeshi MORITA}

\AuthorNameForHeading{T.~Morita}

\Address{Graduate School of Information Science and Technology, Osaka University,
\\
1-1 Machikaneyama-machi, Toyonaka, 560-0043, Japan}
\Email{\href{mailto:t-morita@cr.math.sci.osaka-u.ac.jp}{t-morita@cr.math.sci.osaka-u.ac.jp}}

\ArticleDates{Received October 09, 2012, in f\/inal form July 21, 2013; Published online July 26, 2013}

\Abstract{We show a connection formula for the $q$-conf\/luent hypergeometric functions ${}_2\varphi_1(a,b;0;q,x)$. Combining our connection formula with Zhang's connection formula for ${}_2\varphi_0(a,b;-;q,x)$, we obtain the connection formula for the $q$-conf\/luent hypergeometric equation in the matrix form.  Also we obtain the connection formula of Kummer's conf\/luent hypergeometric functions
 by taking the limit $q\to 1^{-}$ of our connection formula.}

\Keywords{$q$-Borel--Laplace transformation; $q$-dif\/ference equation; connection problem;
$q$-conf\/luent hypergeometric function}

\Classification{33D15; 34M40; 39A13}

\section{Introduction}
We show a~new connection formula for two independent solutions to the $q$-conf\/luent hypergeometric equation
\begin{gather}
(1-abqx)u\big(q^2x\big)-\left\{1-(a+b)qx\right\}u(qx)-qx u(x)=0.
\label{qchge}
\end{gather}

We use notations in accordance with~\cite{GR}.
Assume that $q\in\mathbb{C}^{*}$ satisf\/ies $0<|q|<1$ and $a/b$ is not an integer power of $q$.
The basic hypergeometric series ${}_r\varphi_s$ is def\/ined by
\begin{gather*}
{_r\varphi_s}(a_1,\dots,a_r;b_1,\dots,b_s;q,x):=\sum_{n\ge0}\frac{(a_1,\dots,a_r;q)_n}
{(b_1,\dots,b_s;q)_n(q;q)_n}\left[(-1)^nq^{\frac{n(n-1)}{2}}\right]^{1+s-r}x^n,
\end{gather*}
where  $(a;q)_n$ is the $q$-shifted factorial
\begin{gather*}
(a;q)_n:=
\begin{cases}
1,&
\quad
n=0,
\\
(1-a)(1-aq)\cdots\big(1-aq^{n-1}\big),&
\quad
n\ge1,
\end{cases}
\\
(a;q)_\infty=\lim_{n\to\infty}(a;q)_n,
\end{gather*}
and
\begin{gather*}
(a_1,a_2,\dots,a_m;q)_\infty=(a_1;q)_\infty(a_2;q)_\infty\cdots(a_m;q)_\infty.
\end{gather*}

Equation~\eqref{qchge} has solutions
\begin{gather*}
u_1(x)={}_2\varphi_0(a,b;-;q,x),
\qquad
u_2(x)=\frac{(abx;q)_\infty}{\theta(-qx)}{}_2\varphi_1\left(\frac{q}{a},\frac{q}{b};0;q,abx\right)
\end{gather*}
around the origin and has solutions{\samepage
\begin{gather*}
v_1(x)=x^{-\alpha}{}_2\varphi_1\left(a,0;\frac{aq}{b};q,\frac{q}{abx}\right),
\qquad
v_2(x)=x^{-\beta}{}_2\varphi_1\left(b,0;\frac{bq}{a};q,\frac{q}{abx}\right)
\end{gather*}
around inf\/inity.
Here $q^\alpha=a$ and $q^\beta=b$.}

The connection formula for a~linear $q$-dif\/ference equation of the second order is a~linear relation
between $u_1(x)$, $u_2(x)$ and $v_1(x)$, $v_2(x)$:
\begin{gather*}
\begin{pmatrix}
u_1(x)
\\
u_2(x)
\end{pmatrix}
=
\begin{pmatrix}
C_{11}(x)&C_{12}(x)
\\
C_{21}(x)&C_{22}(x)
\end{pmatrix}
\begin{pmatrix}
v_1(x)
\\
v_2(x)
\end{pmatrix}
,
\end{gather*}
where the connection coef\/f\/icients $C_{jk}(x)$ are $q$-periodic functions.

C.~Zhang~\cite{Z2} proposed the connection formula for $u_1(x)$,
\begin{gather}
{}_2f_0(a,b;\lambda,q,x)=\frac{(b;q)_\infty}{\left(\frac{b}{a};q\right)_\infty}\frac{\theta(a\lambda)}
{\theta(\lambda)}\frac{\theta\left(\frac{qax}{\lambda}\right)}{\theta\left(\frac{qx}{\lambda}\right)}{}
_2\varphi_1\left(a,0;\frac{aq}{b};q,\frac{q}{abx}\right)
\nonumber\\
\phantom{{}_2f_0(a,b;\lambda,q,x)=}
+\frac{(a;q)_\infty}{\left(\frac{a}{b};q\right)_\infty}\frac{\theta(b\lambda)}{\theta(\lambda)}
\frac{\theta\left(\frac{qbx}{\lambda}\right)}{\theta\left(\frac{qx}{\lambda}\right)}{}
_2\varphi_1\left(b,0;\frac{bq}{a};q,\frac{q}{abx}\right),\label{qzhangdiv}
\end{gather}
for $x\in\mathbb{C}^{*}\setminus[-\lambda;q]$.
Here ${}_2f_0(a,b;\lambda,q,x)$ is the $q$-Borel--Laplace transform of the divergent series ${}_2\varphi_0(a,b;-;q,x)$, i.e.,
\[
{}_2f_0(a,b;\lambda ,q,x):=\mathcal{L}_{q,\lambda }^+\circ \mathcal{B}_q^+
{}_2\varphi_0(a,b;-;q,x).
\]

In this paper, we show the following new connection formula for $u_2(x)$:
\begin{gather*}
{}_2\varphi_1(q/a,q/b;0;q,abx)=\frac{(q/a;q)_\infty}{(b/a;q)_\infty}\frac{(aqx,1/ax;q)_\infty}
{(abx;q)_\infty}{}_2\varphi_1(a,0;aq/b;q,1/abx)
\\
\phantom{{}_2\varphi_1(q/a,q/b;0;q,abx)=}
+\frac{(q/b;q)_\infty}{(a/b;q)_\infty}\frac{(bqx,1/bx;q)_\infty}{(abx;q)_\infty}{}_2\varphi_1(b,0;bq/a;q,1/abx).
\end{gather*}

Since $u_1(x)$ is a~divergent series, the $q$-Stokes phenomenon appears in Zhang's connection formula.
But our formula gives the exact relation between the convergent series $u_2(x)$ around the origin and
the convergent series $v_1(x)$, $v_2(x)$ around  inf\/inity.

The theta function of Jacobi is given by the series
\begin{gather*}
\theta_q(x):=\sum_{n\in\mathbb{Z}}q^{\frac{n(n-1)}{2}}x^n,
\qquad
\forall\,  x\in\mathbb{C}^{*},
\end{gather*}
we denote $\theta (x)$ shortly.
The theta function is written by the product form
\[
\theta (x)=\left(q,-x,-\frac{q}{x};q\right)_\infty,
\]
which is known as Jacobi's triple product identity.
For any $k\in\mathbb{Z}$, the theta function also satisf\/ies the $q$-dif\/ference equation
\begin{gather*}
\theta\big(q^kx\big)=q^{-\frac{k(k-1)}{2}}x^{-k}\theta(x).
\end{gather*}
The theta function satisf\/ies the inversion formula $\theta(x)=x\theta(1/x)$.
For all f\/ixed $\lambda\in\mathbb{C}^{*}$, we def\/ine a~$q$-spiral $[\lambda;q]:=\lambda
q^{\mathbb{Z}}=\{\lambda q^k;k\in\mathbb{Z}\}$.
Note that $\theta(\lambda q^k/x)=0$ if and only if $x\in[-\lambda;q]$.

At f\/irst, we review the conf\/luent hypergeometric equation (CHGE).
In 1813~\cite{Gauss13}, C.F.~Gauss studied the hypergeometric series
\begin{gather*}%\label{gauss}
{}_2F_1(\alpha,\beta;\gamma;z)=\sum_{n\ge0}\frac{(\alpha)_n(\beta)_n}{(\gamma)_nn!}z^n,
\qquad
\gamma\not=0,-1,-2,\dots,
\end{gather*}
where $(\alpha)_n=\alpha\{\alpha+1\}\cdots\{\alpha+(n-1)\}$.

More generally, the generalized hypergeometric series is given by
\begin{gather*}
{}_rF_s(\alpha_1,\dots,\alpha_r;\beta_1,\dots,\beta_s;z)=\sum_{n\ge0}\frac{(\alpha_1)_n\cdots(\alpha_r)_n}
{(\beta_1)_n\cdots(\beta_s)_nn!}z^n.
\end{gather*}
The hypergeometric function ${}_2F_1(\alpha,\beta;\gamma;z)$ satisf\/ies the second-order dif\/ferential equation
\begin{gather}
\label{gausseq}
z(1-z)\frac{d^2u}{dz^2}+\left\{\gamma-(\alpha+\beta+1)z\right\}\frac{du}{dz}-\alpha\beta u=0.
\end{gather}
Gauss gave the connection formula for the function ${}_2F_1(\alpha,\beta;\gamma,z)$.
We put $z\mapsto z/\beta$, take the limit $\beta\to\infty$ in equation~\eqref{gausseq}, and obtain
the conf\/luent hypergeometric equation (CHGE)
\begin{gather}
\label{CHGE}
z\frac{d^2u}{dz^2}+(\gamma-z)\frac{du}{dz}-\alpha u=0.
\end{gather}
Solutions of~\eqref{CHGE} around the origin are
\begin{gather*}%\label{solo01}
\hat{u}_1(z)={}_1F_1(\alpha;\gamma;z)
\end{gather*}
and
\begin{gather}
\hat{u}_2(z)=z^{1-\gamma}{}_1F_1(\alpha-\gamma+1,2-\gamma,z).
\label{solo02}
\end{gather}
Solutions around inf\/inity are given by the divergent series
\begin{gather*}
\hat{v}_1(z)=(-z)^{-\alpha}{}_2F_0(\alpha,\alpha-\gamma+1;-;1/z)
\end{gather*}
and
\begin{gather*}
\hat{v}_2(z)=e^zz^{\alpha-\gamma}{}_2F_0(1-\alpha,\gamma-\alpha;-;1/z).
\end{gather*}
The asymptotic expansion of ${}_1F_1(\alpha;\gamma;z)$ is given by
\begin{gather}
{}_1F_1(\alpha;\gamma;z)\sim\frac{\Gamma(\gamma)}{\Gamma(\gamma-\alpha)}(-z)^{-\alpha}{}
_2F_0(\alpha,\alpha-\gamma+1;-;-1/z)
\nonumber\\
\phantom{{}_1F_1(\alpha;\gamma;z)=}
+\frac{\Gamma(\gamma)}{\Gamma(\alpha)}e^zz^{\alpha-\gamma}{}_2F_0(1-\alpha,\gamma-\alpha;-;1/z),\label{eq6}
\end{gather}
where $-\pi/2<\arg z<3\pi/2$.
Note that the connection formula for the second solution around inf\/inity~\eqref{solo02} can be
derived from \eqref{eq6}.
In Section~\ref{Section2} we deal with another degeneration of equation~\eqref{gausseq} which is slightly dif\/ferent
from the standard way.

It is known that there exists a $q$-analogue of ${}_2F_1(\alpha,\beta;\gamma;z)$, which was introduced by
E.~Heine in 1847 as
\begin{gather*}
{}_2\varphi_1(a,b;c;q,x):=\sum_{n\ge0}\frac{(a;q)_n(b;q)_n}{(c;q)_n(q;q)_n}x^n.
\end{gather*}
We assume that $c$ is not integer powers of $q$.
The function ${}_2\varphi_1(a,b;c;q,x)$ satisf\/ies the second-order $q$-dif\/ference equation ($q$-HGE)
\begin{gather}
\label{qgauss}
x(c-abqx)\mathcal{D}_q^2u+\left[\frac{1-c}{1-q}+\frac{(1-a)(1-b)-(1-abq)}{1-q}x\right]\mathcal{D}_qu
-\frac{(1-a)(1-b)}{(1-q)^2}u=0,
\end{gather}
where $\mathcal{D}_q$ is the $q$-derivative operator def\/ined for f\/ixed $q$ by
\begin{gather*}
\mathcal{D}_qf(x)=\frac{f(x)-f(qx)}{(1-q)x}.
\end{gather*}
By replacing $a$, $b$, $c$ by $q^\alpha$, $q^\beta$, $q^\gamma$ and then letting $q\to1^-$,
equation~\eqref{qgauss} tends to the hypergeometric equation~\eqref{gausseq}.
The $q$-hypergeometric equation~\eqref{qgauss} can be rewritten as
\begin{gather}
\label{qqgauss}
(c-abqx)u\big(q^2x\big)-\left\{c+q-(a+b)qx\right\}u(qx)+q(1-x)u(x)=0.
\end{gather}
If we set $x\mapsto cx$ and $c\to\infty$ in~\eqref{qqgauss}, we obtain the $q$-conf\/luent hypergeometric
equation~\eqref{qchge}.
Equation~\eqref{qchge} is considered as a~$q$-analogue of CHGE.
Note that the f\/irst solution $u_1(x)$ is a~divergent series and $u_2(x)$ is a~convergent series
around the origin.
Therefore we should study the connection formula for~$u_1(x)$ and~$u_2(x)$ independently.
We need dif\/ferent types of $q$-Borel--Laplace transformations to
obtain the connection formula for~$u_1(x)$ and~$u_2(x)$.
This point is essentially dif\/ferent from the dif\/ferential equation case.

We study connection problems for
linear $q$-dif\/ference equations with irregular singular points.
The irregularity of $q$-dif\/ference equations are studied using the Newton polygons by J.-P.~Ramis, J.~Sauloy and C.~Zhang~\cite{RSZ}. For any $q$-dif\/ference operator $\mathcal{P}=\sigma_q^n+a_1(z)\sigma_q^{n-1}+\cdots +a_n(z)$, the Newton polygon is def\/ined as the convex hull of $\{(i,j)\in\mathbb{Z}^2\,|\,j\ge v_0(a_i)\}$, provided that $v_0$ are $z$-adic valuation in suitable fields. Graphically, the irregularity of $q$-dif\/ference equations and $q$-dif\/ference modules are the height of the Newton polygon (from the bottom to the upper right end).

Connection problems for
linear $q$-dif\/ference equations~\cite{Ohyama} with regular
singular points were studied by G.D.~Birkhof\/f~\cite{Birkhoff}.
Linear $q$-dif\/ference equations have formal power series solutions $x^{\alpha} \sum\limits_{n\ge0}a_nx^n$
around the origin and $x^\beta \sum\limits_{n\ge0}b_n x^{-n}$ around  inf\/inity for generic exponents.
But for
connection problems for
linear $q$-dif\/ference equations, we replace the function~$x^\kappa$ with
the function $\theta(x)/\theta(k x)$, where $k=q^\kappa$,
since these functions satisfy the same
$q$-dif\/ference equation $\sigma_qf(x)=q^{\kappa}f(x)$.
Then, the fundamental system of solutions is given by single valued functions which have single poles at
suitable $q$-spirals.
Therefore, each connection coef\/f\/icient has the periods $q$ and $e^{2\pi i}$.
The f\/irst example was given by G.N.~Watson~\cite{W}.
But a~few examples of irregular singular cases are known~\cite{M0, Z0,Z1}.
In this paper, we give a~connection formula for the $q$-conf\/luent type function using the
$q$-Borel--Laplace transformations.

In 1910, Watson~\cite{W} showed that the connection formula for the series ${}_2\varphi_1(a,b;c;q,x)$ has the following form
\begin{gather}
{}_2\varphi_1\left(a,b;c;q,x\right)=\frac{(b,c/a;q)_\infty}{(c,b/a;q)_\infty}\frac{\theta(-ax)}{\theta(-x)}
{}_2\varphi_1\left(a,\frac{aq}{c};\frac{aq}{b};q,\frac{cq}{abx}\right)
\nonumber
\\
\phantom{{}_2\varphi_1\left(a,b;c;q;x\right)=}
+\frac{(a,c/b;q)_\infty}{(c,a/b;q)_\infty}\frac{\theta(-bx)}{\theta(-x)}{}_2\varphi_1\left(b,\frac{bq}{c}
;\frac{bq}{a};q,\frac{cq}{abx}\right).
\label{watson's2phi1}
\end{gather}
Note that we can not set $a=0$ or $b=0$ directly in this formula.

In 2002, C.~Zhang~\cite{Z2} showed one of the connection formula for the $q$-CHGE \eqref{qzhangdiv}.
The $q$-Borel--Laplace transformations were studied by C.~Zhang in~\cite{Z2} (see Section~\ref{Section2} for more details).
When we study connection problems for
$q$-dif\/ference equations, this resummation method becomes a~powerful tool.
Note that we can f\/ind a~new parameter $\lambda$ in the resummation ${}_2f_0(a,b;\lambda,q,x)$.
Here $\lambda$ is the direction of the summation.
This parameter brings us new viewpoints for the study of the $q$-Stokes phenomenon.

It is known that there exist two dif\/ferent types of the $q$-Borel--Laplace transformations.

The $q$-Borel--Laplace transformations of the f\/irst kind are def\/ined in~\cite{Z2} and the
$q$-Borel--Laplace transformations of the second kind are studied in~\cite{Z1}.
These $q$-Borel transformations are formal inverse transformations of each of the $q$-Laplace
transformations.

C.~Zhang presented a~connection formula for the series ${}_2\varphi_0(a,b;-;q,x)$ by the $q$-Borel--Laplace
transformations of the f\/irst kind.
But the connection formula for the second solution of~\eqref{qchge} is not known.
In this paper, we show the second connection formula for $q$-CHGE with the using of the $q$-Borel
transformation and the $q$-Laplace transformation of the second kind.
Combining with Zhang's connection formula, we obtain the connection formula in the matrix form (see Theorem~\ref{theorem_main}).
Using Watson's formula~\eqref{watson's2phi1} we also give another proof of the new connection formula in Section~\ref{ano}.

In Section~\ref{Section3} we consider the limit $q\to1^-$ of our connection formula.
If we take the limit $q\to1^-$, we formally obtain the connection formula for the
conf\/luent hypergeometric series ${}_2F_0$.

\section{The connection formula and the connection matrix}\label{Section2}

We review a~$q$-conf\/luent hypergeometric
equation in Section~\ref{Subsection21}.
Then we show a~connection formula for the $q$-conf\/luent hypergeometric function, which is dif\/ferent
from Zhang's formula.

\subsection{Conf\/luent hypergeometric equation}\label{Subsection21}

For the conf\/luent hypergeometric
equation~\eqref{gausseq}, we take another degeneration.
We put $z\mapsto z\gamma$ and take the limit $\gamma\to\infty$.
Then we obtain
\begin{gather}
\label{difcon}
z^2\frac{d^2u}{dz^2}-\left\{1-(\alpha+\beta+1)z\right\}\frac{du}{dz}+\alpha\beta u=0.
\end{gather}
Solutions to~\eqref{difcon} around the origin are given by the divergent series
\begin{gather*}
\tilde{u}_1(z)={}_2F_0(\alpha,\beta;-,z)
\qquad \mbox{and} \qquad
\tilde{u}_2(z)=e^{\frac{1}{z}}(-z)^{1-\alpha-\beta}\, {}_2F_0(1-\alpha,1-\beta;-,z).
\end{gather*}
Solutions around inf\/inity are given by the convergent series
\begin{gather*}
\tilde{v}_1(z)=(-z)^\alpha\, {}_1F_1(\alpha,1+\alpha-\beta,-1/z)
\qquad\! \mbox{and}\!
\qquad
\tilde{v}_2(z)=(-z)^\beta\, {}_1F_1(\beta,1+\beta-\alpha,-1/z).
\end{gather*}

We consider a $q$-analogue of the conf\/luent hypergeometric equation~\eqref{qqgauss}.  The second-order $q$-dif\/ference equation
\begin{gather}
x\left\{abqx-(1-q)\right\}D_q^2u(x) +\left\{1-\frac{(1-a)(1-b)-(1-abq)}{1-q}x\right\}D_qu(x)\nonumber\\
 \qquad{} +\frac{(1-a)(1-b)}{(1-q)^2}u(x)=0. \label{eq9}
\end{gather}
can be rewritten as
\begin{gather}\label{con1}
(1-abqx)u\big(xq^2\big)-\left\{1-(a+b)qx\right\}u(xq)-qxu(x)=0,
\end{gather}
which is called a~$q$-conf\/luent hypergeometric equation. When we take $q\to 1^-$, the limit of \eqref{eq9} is the dif\/ferential equation~\eqref{qqgauss}, provided that $a=q^\alpha$, $b=q^\beta$.

\subsection[Local solutions to the $q$-conf\/luent hypergeometric equation]{Local solutions to the $\boldsymbol{q}$-conf\/luent hypergeometric equation}

Consider the connection problem of~\eqref{con1}.
At f\/irst we show local solutions for~\eqref{con1} around $x=0$ and $x=\infty$.
\begin{lemma}
Equation~\eqref{con1} has  solutions
\begin{gather}
u_1(x)={}_2\varphi_0(a,b;-;q,x),
\label{solo1}
\\
u_2(x)=\frac{(abx;q)_\infty}{\theta(-qx)}{}_2\varphi_1\left(\frac{q}{a},\frac{q}{b};0;q,abx\right)
\label{solo2}
\end{gather}
around the origin and  solutions
\begin{gather}
v_1(x)=x^{-\alpha}{}_2\varphi_1\left(a,0;\frac{aq}{b};q,\frac{q}{abx}\right),
\label{soli1}
\\
v_2(x)=x^{-\beta}{}_2\varphi_1\left(b,0;\frac{bq}{a};q,\frac{q}{abx}\right),
\label{soli2}
\end{gather}
around   inf\/inity, provided that $a=q^\alpha$ and $b=q^\beta$.
\end{lemma}
\begin{proof}
We show a~fundamental system of solutions of~\eqref{con1} around $x=0$.
If set $u(x)=\sum\limits_{n\ge0}a_n x^n$, $a_0=1$, then we obtain
\begin{gather*}
u_1(x)={}_2\varphi_0(a,b;-;q,x).
\end{gather*}
We set $\mathcal{E}(x)=1/\theta(-qx)$ and $f(x)=\sum\limits_{n\ge0}a_n x^n$, $a_0=1$ to obtain another solution solution
around the origin.
We assume that $u(x)=\mathcal{E}(x)f(x)$.
Note that the function $\mathcal{E}(x)$ has the following property
\begin{gather*}
\sigma_q\mathcal{E}(x)=-qx\mathcal{E}(x),
\qquad
\sigma_q^2\mathcal{E}(x)=q^3x^2\mathcal{E}(x).
\end{gather*}
Therefore, we obtain the equation
\begin{gather}
\left[q^3x(1-abqx)\sigma_q^2+q\left\{1-(a+b)qx\right\}\sigma_q-q\right]f(x)=0.
\label{con2}
\end{gather}
Since the inf\/inite product $(abx;q)_\infty$ satisf\/ies the following $q$-dif\/ference relation
\begin{gather*}
\sigma_q\left[(abx;q)_\infty\right]=\frac{1}{1-abx}(abx;q)_\infty,
\end{gather*}
we obtain the second solution.
Therefore, solutions of equation~\eqref{con1} around the origin are given~by
\begin{gather*}
u_1(x)={}_2\varphi_0(a,b;-;q,x),
\qquad
u_2(x)=\frac{(abx;q)_\infty}{\theta(-qx)}{}_2\varphi_1\left(\frac{q}{a},\frac{q}{b};0;q,abx\right).
\end{gather*}

Around $x=\infty$, we can easily determine local solutions by setting
\begin{gather*}
v(x)=\frac{\theta(a\mu x)}{\theta(\mu x)}\sum_{n\ge0}a_n x^{-n},
\qquad
a_0=1,
\end{gather*}
for any f\/ixed $\mu\in\mathbb{C}^{*}$ and $x\in\mathbb{C}^{*}\setminus[-\mu;q]$.
\end{proof}

Here $u_1(x)$ is a~divergent series and $u_2(x)$, $v_1(x)$ and $v_2(x)$ are convergent series~\cite{GR}.
Therefore, the $q$-Stokes phenomenon occurs for $u_1(x)$.
\begin{definition}
For any $f(x)=\sum\limits_{n\ge0}a_nx^n$, the $q$-Borel transformation $\mathcal{B}_q^+$ is
\begin{gather*}
\left(\mathcal{B}_q^+f\right)(\xi)=\varphi(\xi):=\sum_{n\ge0}a_nq^{\frac{n(n-1)}{2}}\xi^n,
\end{gather*}
and the $q$-Laplace transformation $\mathcal{L}_{q,\lambda}^+$ is
\begin{gather*}
\big(\mathcal{L}_{q,\lambda}^+\varphi\big)(x):=\sum_{n\in\mathbb{Z}}\frac{\varphi(q^n\lambda)}
{\theta\left(\frac{q^n\lambda}{x}\right)}.
\end{gather*}
\end{definition}
C.~Zhang determined a~resummation of~\eqref{solo1} by the $q$-Borel--Laplace transformations of the
f\/irst kind as follows
\begin{gather*}
{}_2f_0(a,b;\lambda,q,x):=\mathcal{L}_{q,\lambda}^+\circ\mathcal{B}_q^+{}_2\varphi_0(a,b;-;q,x).
\end{gather*}
He also presented a~connection formula~\eqref{qzhangdiv} for ${}_2f_0(a,b;\lambda,q,x)$.

But the connection formula between~\eqref{solo2} and~\eqref{soli1}, \eqref{soli2} is not known.
In the next section, we show the second connection formula by means of the $q$-Borel--Laplace
transformations of the second kind.

\subsection{The second connection formula}

We def\/ine the $q$-Borel transformation and the $q$-Laplace transformation
of the second kind. These transformations are introduced by C.~Zhang to obtain the solution of equation~\eqref{con2}.
\begin{definition}
For $f(x)=\sum\limits_{n\ge0}a_nx^n$, the $q$-Borel transformation is def\/ined by
\begin{gather*}
g(\xi)=\left(\mathcal{B}^-_qf\right)(\xi):=\sum_{n\ge0}a_nq^{-\frac{n(n-1)}{2}}\xi^n,
\end{gather*}
and the $q$-Laplace transformation is
\begin{gather*}
\left(\mathcal{L}_q^-g\right)(x):=\frac{1}{2\pi i}\int_{|\xi|=r}g(\xi)\theta_q\left(\frac{x}{\xi}
\right)\frac{d\xi}{\xi}.
\end{gather*}
\end{definition}
Here $r>0$ is a suf\/f\/iciently small number.
The $q$-Borel transformation is considered as a~formal inverse of the $q$-Laplace transformation.
\begin{lemma}[\cite{Z1}]
We assume that the function $f$ can be $q$-Borel transformed to the analytic function
$g(\xi)$ around $\xi=0$.
Then, we have
\begin{gather*}
\mathcal{L}^-_q\circ\mathcal{B}^-_qf=f.
\end{gather*}
\end{lemma}
\begin{proof}
We can prove this lemma calculating residues of the $q$-Laplace transformation around the origin.
\end{proof}
The $q$-Borel transformation satisf\/ies the following operational relation.
\begin{lemma}
\label{oprel}
For any $l,m\in\mathbb{Z}_{\ge0}$,
\begin{gather*}
\mathcal{B}^-_q\big(x^m\sigma_q^l\big)=q^{-\frac{m(m-1)}{2}}\xi^m\sigma_q^{l-m}\mathcal{B}^-_q.
\end{gather*}
\end{lemma}
We apply the $q$-Borel transformation to equation~\eqref{con2} and use Lemma~\ref{oprel}.
We use the notation~$g(\xi)$ as the $q$-Borel transform of $u_2(x)$.
We check out that $g(\xi)$ satisf\/ies the f\/irst-order $q$-dif\/ference equation
\begin{gather*}
g(q\xi)=\frac{(1+aq\xi)(1+bq\xi)}{(1+q^2\xi)}g(\xi).
\end{gather*}
Since $g(0)=a_0=1$, we have the inf\/inite product of $g(\xi)$ as follows
\begin{gather*}
g(\xi)=\frac{(-q^2\xi;q)_\infty}{(-qa\xi;q)_\infty(-qb\xi;q)_\infty}.
\end{gather*}
Note that $g(\xi)$ has single poles at
\begin{gather*}
\left\{\xi\in\mathbb{C}^{*};\xi=-\frac{1}{aq^{k+1}},-\frac{1}{bq^{k+1}},k\in\mathbb{Z}_{\ge0}\right\}.
\end{gather*}
We set
\begin{gather*}
r_0:=\max\left\{\frac{1}{|aq|},
\frac{1}{|bq|}\right\}
\end{gather*}
and choose the radius $r>0$ such that $0<r<r_0$.
By Cauchy's residue theorem, the $q$-Laplace transform of $g(\xi)$ is
\begin{gather*}
f(x)
=\frac{1}{2\pi i}\int_{|\xi|=r}g(\xi)\theta\left(\frac{x}{\xi}\right)\frac{d\xi}{\xi}
\\
\phantom{f(x)}
=-\sum_{k\ge0}\operatorname{Res}\left\{g(\xi)\theta\left(\frac{x}{\xi}\right)\frac{1}{\xi};\xi=-\frac{1}
{aq^{k+1}}\right\}
-\sum_{k\ge0}\operatorname{Res}\left\{g(\xi)\theta\left(\frac{x}{\xi}\right)\frac{1}{\xi};\xi=-\frac{1}
{bq^{k+1}}\right\},
\end{gather*}
where $0<r<r_0$.
Since there exists a~positive constant $C_N$ (for any integer $N$) s.t.,
\begin{gather*}
|g(\xi)|\le C_N\xi^{-N}.
\end{gather*}
The following lemma plays a key role to calculate the residue.
\begin{lemma}
For any $k\in\mathbb{N}$, $\lambda\in\mathbb{C}^{*}$, we have:
\begin{gather*}
1) \ \ \operatorname{Res}\left\{\frac{1}{\left(\xi/\lambda;q\right)_\infty}\frac{1}{\xi}:\xi=\lambda
q^{-k}\right\}=\frac{(-1)^{k+1}q^{\frac{k(k+1)}{2}}}{(q;q)_k(q;q)_\infty} ,\\
2) \ \   \frac{1}{(\lambda
q^{-k};q)_\infty}=\frac{(-\lambda)^{-k}q^{\frac{k(k+1)}{2}}}{(\lambda;q)_\infty\left(q/\lambda;q\right)_k},
\qquad
\lambda\not\in q^{\mathbb{Z}} .
\end{gather*}
\end{lemma}
Summing up all
residues, we obtain $f(x)$ as follows
\begin{gather*}
f(x)= \frac{\left(\frac{q}{a};q\right)_\infty}{\left(\frac{b}{a},q;q\right)_\infty}\frac{\theta(-aqx)}
{\theta(-qx)}{}_2\varphi_1\!\!\left(a,0;\frac{aq}{b};q,\frac{q}{abx}\right)
+\frac{\left(\frac{q}{b};q\right)_\infty}{\left(\frac{a}{b},q;q\right)_\infty}\frac{\theta(-bqx)}
{\theta(-qx)}{}_2\varphi_1\!\!\left(b,0;\frac{bq}{a};q,\frac{q}{abx}\right)\!.
\end{gather*}
Therefore, we obtain the following theorem.
\begin{theorem}
\label{thm1}
For any $x\not\in[1;q]$, we have
\begin{gather}
u_2(x)=\frac{(abx;q)_\infty}{\theta(-qx)}{}_2\varphi_1\left(\frac{q}{a},\frac{q}{b};0;q,abx\right)
=\frac{\left(\frac{q}{a};q\right)_\infty}{\left(\frac{b}{a},q;q\right)_\infty}\frac{\theta(-aqx)}
{\theta(-qx)}{}_2\varphi_1\left(a,0;\frac{aq}{b};q,\frac{q}{abx}\right)\nonumber\\
\phantom{u_2(x)=}
{}+\frac{\left(\frac{q}{b};q\right)_\infty}{\left(\frac{a}{b},q;q\right)_\infty}\frac{\theta(-bqx)}
{\theta(-qx)}{}_2\varphi_1\left(b,0;\frac{bq}{a};q,\frac{q}{abx}\right).\label{eq16}
\end{gather}
\end{theorem}

\subsection{The connection matrix}

Combining Zhang's connection formula and Theorem~\ref{thm1}, we give the connection matrix for equation~\eqref{qchge}. At f\/irst, we def\/ine a new fundamental system of solutions around inf\/inity.
For any $\lambda$, $\mu\in\mathbb{C}^{*}$, $x\in\mathbb{C}^{*}\setminus[-\mu;q]$, $S_\mu(a,b;q,x)$ is
\begin{gather*}
S_\mu(a,b;q,x):=\frac{\theta(a\mu x)}{\theta(\mu x)}{}_2\varphi_1\left(a,0;\frac{aq}{b};q,\frac{q}{abx}
\right).
\end{gather*}
The function $\theta(a\mu x)/\theta(\mu x)$ satisf\/ies the following $q$-dif\/ference equation
\begin{gather*}
u(qx)=\frac{1}{a}u(x),
\end{gather*}
which is also satisf\/ied by the function $u(x)=x^{-\alpha}$, $a=q^{\alpha}$.
Note that the pair $\left(S_\mu(a,b;q,x)\right.$, $\left.S_\mu(b,a;q,x)\right)$ gives a~fundamental system of
solutions of equation~\eqref{qchge} if $a/b\not\in q^{\mathbb{Z}}$.
We def\/ine $q$-elliptic functions $C_\mu^{\lambda}(a,b;q,x)$ and $C_\mu (a,b;q,x)$.

\begin{definition}
For any $\lambda,\mu\in\mathbb{C}^{*}$
we set functions $C_\mu^\lambda(a,b;q,x)$ and $C_\mu(a,b;q,x)$ as follows
\begin{gather*}
C_\mu^\lambda(a,b;q,x):=\frac{\left(b;q\right)_\infty}{\left(\frac{b}{a};q\right)_\infty}
\frac{\theta(a\lambda)}{\theta(\lambda)}\frac{\theta\left(\frac{qax}{\lambda}\right)}{\theta\left(\frac{qx}
{\lambda}\right)}\frac{\theta(\mu x)}{\theta(a\mu x)},
\\
C_\mu(a,b;q,x):=\frac{\left(\frac{q}{a};q\right)_\infty}{\left(\frac{b}{a},q;q_{\infty}\right)}
\frac{\theta(-aqx)}{\theta(-qx)}\frac{\theta(\mu x)}{\theta(a\mu x)}.
\end{gather*}
\end{definition}

Then $C_\mu^\lambda(a,b;q,x)$ and $C_\mu(a,b;q,x)$ are single valued as functions of~$x$.
They satisfy the follo\-wing relation
\begin{gather*}
C_\mu^\lambda(a,b;q,e^{2\pi i}x)=C_\mu^\lambda(a,b;q,x),
\qquad
C_\mu^\lambda(a,b;q,qx)=C_\mu^\lambda(a,b;q,x)
\end{gather*}
and
\begin{gather*}
C_\mu\big(a,b;q,e^{2\pi i}x\big)=C_\mu(a,b;q,x),
\qquad
C_\mu(a,b;q,qx)=C_\mu(a,b;q,x).
\end{gather*}
i.e.\ $C_\mu^\lambda(a,b;q,x)$ and $C_\mu(a,b;q,x)$ are $q$-elliptic functions.
We set
\begin{gather*}
{}_2f_1(a,b;q,x):=u_2(x)=\frac{(abx;q)_\infty}{\theta(-qx)}{}_2\varphi_1\left(\frac{q}{a},\frac{q}{b}
;0;q,abx\right).
\end{gather*}
Thus, we obtain the connection formula in the matrix form.
\begin{theorem}\label{theorem_main}
For any $\lambda,\mu\in\mathbb{C}^{*}$, $x\in\mathbb{C}^{*}\setminus[1;q]\cup[-\mu/a;q]\cup[-\mu/b;q]\cup[-\lambda;q]$, we have
\begin{gather*}
\begin{pmatrix}
{}_2f_0(a,b;\lambda,q,x)
\\
{}_2f_1(a,b;q,x)
\end{pmatrix}
=
\begin{pmatrix}
C_\mu^\lambda(a,b;q,x)&C_\mu^\lambda(b,a;q,x)
\\
C_\mu(a,b;q,x)&C_\mu(b,a;q,x)
\\
\end{pmatrix}
\begin{pmatrix}
S_\mu(a,b;q,x)
\\
S_\mu(b,a;q,x)
\end{pmatrix}
.
\end{gather*}
\end{theorem}

\subsection{Derivation from Watson's formula}
\label{ano}
In this section, we give another proof of Theorem~\ref{thm1}.
Watson's formula~\eqref{watson's2phi1} is a connection formula for the basic hypergeometric functions ${}_2\varphi_1(a,b;c;q,x)$.
We derive the connection formula
\begin{gather*}
\frac{(abx;q)_\infty}{\theta(-qx)}{}_2\varphi_1\left(\frac{q}{a},\frac{q}{b};0;q,abx\right)
=\frac{\left(q/a;q\right)_\infty}{\left(b/a,q;q\right)_\infty}\frac{\theta(-aqx)}{\theta(-qx)}\, {}
_2\varphi_1\left(a,0;\frac{aq}{b};q,\frac{q}{abx}\right)
\\
\hphantom{\frac{(abx;q)_\infty}{\theta(-qx)}{}_2\varphi_1\left(\frac{q}{a},\frac{q}{b};0;q,abx\right)=}{}
+\frac{\left(q/b;q\right)_\infty}{\left(a/b,q;q\right)_\infty}\frac{\theta(-bqx)}{\theta(-qx)}\, {}
_2\varphi_1\left(b,0;\frac{bq}{a};q,\frac{q}{abx}\right)
\end{gather*}
from Watson's formula \eqref{watson's2phi1}.
By take the limit $c\to 0$ of Watson's formula, we obtain the following proposition.
\begin{proposition}[\cite{Z2}]
For any $x\in\mathbb{C}^{*}\setminus q^{\mathbb{Z}}$, we have
\begin{gather}
{}_2\varphi_1(a,b;0;q,x)=
\frac{(b;q)_\infty}{(b/a;q)_\infty}\frac{\theta(-ax)}{\theta(-x)}{}
_1\varphi_1\left(a;\frac{aq}{b};q,\frac{q^2}{bx}\right)
\nonumber
\\
\phantom{{}_2\varphi_1(a,b;0;q,x)=}
+\frac{(a;q)_\infty}{(a/b;q)_\infty}\frac{\theta(-bx)}{\theta(-x)}{}_1\varphi_1\left(b;\frac{bq}{a};q,\frac{q^2}{ax}\right).
\label{wat2}
\end{gather}
\end{proposition}
In~\eqref{wat2}, we put $a\mapsto q/a$, $b\mapsto q/b$ and $x\mapsto abx$, we obtain the relation as
follows
\begin{corollary}
For any $x\in\mathbb{C}^{*}\setminus[ab;q]$, we have
\begin{gather}
{}_2\varphi_1\left(\frac{q}{a},\frac{q}{b};0;q,abx\right)
=\frac{(q/b;q)_\infty}{(a/b;q)_\infty}
\frac{\theta(-bqx)}{\theta(-abx)}{}_1\varphi_1\left(\frac{q}{a};\frac{bq}{a};q,\frac{q}{ax}\right)
\nonumber
\\
\phantom{{}_2\varphi_1\left(\frac{q}{a},\frac{q}{b};0;q,abx\right)=}
+\frac{(q/a;q)_\infty}{(b/a;q)_\infty}\frac{\theta(-aqx)}{\theta(-abx)}{}_1\varphi_1\left(\frac{q}{b}
;\frac{aq}{b};q,\frac{q}{bx}\right).
\label{wat3}
\end{gather}
\end{corollary}
The function ${}_1\varphi_1(a;c;q,x)$ is related to~${}_2\varphi_1(a',0;c';q,x)$.
Consider the relation between the function~${}_1\varphi_1$ and the function~${}_2\varphi_1$.
\begin{proposition}
For any $x\in\mathbb{C}^{*}$, we have
\begin{gather}
\label{prpr}
{}_1\varphi_1\left(a_1;c_1;q,\frac{c_1x}{a_1}\right)=(x;q)_\infty{}_2\varphi_1\left(\frac{c_1}{a_1}
,0;c_1;q,x\right),
\end{gather}
provided that $c_1/a_1\not\in q^{\mathbb{Z}}$.
\end{proposition}
\begin{proof}
The function ${}_1\varphi_1(a_1;c_1;q,c_1x/a_1)$ satisf\/ies the equation
\begin{gather*}
\left[(c_1-c_1qx)\sigma_q^2-\left\{(c_1+q)-\frac{qc_1}{a_1}x\right\}\sigma_q+q\right]v(x)=0.
\end{gather*}
We set $v(x)=(x;q)_\infty\tilde{v}(x)$, where $\tilde{v}(x):=\sum\limits_{n\ge0}\tilde{v}_nx^n$ and
$\tilde{v}_0:=1$.
Note that the function $(x;q)_\infty$ satisf\/ies the f\/irst-order $q$-dif\/ference equation
\begin{gather*}
\sigma_qf(x)=\frac{1}{1-x}f(x).
\end{gather*}
Then, we obtain the equation
\begin{gather}
\label{vofx}
\left[\sigma_q^2-\left\{\left(1+\frac{q}{c_1}\right)-\frac{qx}{a_1}\right\}\sigma_q+\frac{q}{c_1}
(1-x)\right]\tilde{v}(x)=0.
\end{gather}
Equation~\eqref{vofx} has the solution
\begin{gather*}
\tilde{v}(x)={}_2\varphi_1\left(\frac{c_1}{a_1},0;c_1;q,x\right).
\end{gather*}
Therefore, we obtain the conclusion.
\end{proof}
\begin{corollary}
In~\eqref{prpr}, we put $a_1\mapsto q/a$, $c_1\mapsto bq/a$ and $x\mapsto q/abx$. Then we obtain
\begin{gather}
\label{c1}
{}_1\varphi_1\left(\frac{q}{a};\frac{bq}{a};q,\frac{q}{ax}\right)=\left(\frac{q}{abx};q\right)_\infty{}
_2\varphi_1\left(b,0;\frac{bq}{a};q,\frac{q}{abx}\right).
\end{gather}
We also put $a_1\mapsto q/b$, $c_1\mapsto aq/b$ and $x\mapsto q/abx$. Then we obtain
\begin{gather}
\label{c2}
{}_1\varphi_1\left(\frac{q}{b};\frac{aq}{b};q,\frac{q}{bx}\right)=\left(\frac{q}{abx};q\right)_\infty{}
_2\varphi_1\left(a,0;\frac{aq}{b};q,\frac{q}{abx}\right).
\end{gather}
\end{corollary}
By relations~\eqref{wat3},~\eqref{c1} and~\eqref{c2},
\begin{gather*}
\frac{(abx;q)_\infty}{\theta(-qx)}{}_2\varphi_1\left(\frac{q}{a},\frac{q}{b};0;q,abx\right)\\
\qquad{}
=\frac{(q/b;q)_\infty}{(a/b,q;q)_\infty}\frac{\theta(-bqx)}{\theta(-abx)}\frac{\left(q,abx,\frac{q}{abx}
;q\right)_\infty}{\theta(-qx)}\, {}_2\varphi_1\left(b,0;\frac{bq}{a};q,\frac{q}{abx}\right)
\\
\qquad\quad{}
+\frac{(q/a;q)_\infty}{(b/a,q;q)_\infty}\frac{\theta(-aqx)}{\theta(-abx)}\frac{\left(q,abx,\frac{q}{abx}
;q\right)_\infty}{\theta(-qx)}\, {}_2\varphi_1\left(a,0;\frac{aq}{b};q,\frac{q}{abx}\right)
\\
\qquad{}
=\frac{(q/b;q)_\infty}{(a/b,q;q)_\infty}\frac{\theta(-bqx)}{\theta(-qx)}\, {}_2\varphi_1\left(b,0;\frac{bq}{a}
;q,\frac{q}{abx}\right)
\\
\qquad\quad{}
+\frac{(q/a;q)_\infty}{(b/a,q;q)_\infty}\frac{\theta(-aqx)}{\theta(-qx)}\, {}_2\varphi_1\left(a,0;\frac{aq}{b}
;q,\frac{q}{abx}\right).
\end{gather*}
Therefore, we obtain the formula~\eqref{eq16}.

\section[The limit $q\to1^-$ of the connection formula]{The limit $\boldsymbol{q\to1^-}$ of the connection formula}\label{Section3}

In this section we show the limit $q\to1^-$ of our connection formula.
In~\cite{Z2}, C.~Zhang proposed the following limit.
\begin{theorem}[\cite{Z2}]%\label{Zhang}
For any $\alpha,\beta\in\mathbb{C}^{*}(\alpha-\beta\not\in\mathbb{Z})$ and $z$ in any compact domain of
$\mathbb{C}^{*}\setminus[-\infty,0]$, we have
\begin{gather*}
\lim_{q\to1^-}{}_2f_0\left(q^\alpha,q^\beta;\lambda,q,\frac{z}{(1-q)}\right)
\\
\qquad
=\frac{\Gamma(\beta-\alpha)}
{\Gamma(\beta)}z^{-\alpha}{}_1F_1\left(\alpha,\alpha-\beta+1;\frac{1}{z}\right)
+\frac{\Gamma(\alpha-\beta)}{\Gamma(\alpha)}z^{-\beta}{}_1F_1\left(\beta;\beta-\alpha+1;\frac{1}{z}\right).
\end{gather*}
\end{theorem}

Our limit of the connection formula in Theorem~\ref{thm1} dif\/fers from the theorem above.
By Theorem~\ref{thm1}, we have{\samepage
\begin{gather}
\label{moth1}
u_2(x)=C_\mu(a,b;q,x)S_\mu(a,b;q,x)+C_\mu(b,a;q,x)S_\mu(b,a;q,x)
\end{gather}
for any $(x,q)\in\mathbb{C}^{*}\times(0,1]$.}

The limit $q\to1^-$ of the left-hand side of~\eqref{moth1} is formally given by
$e^{1/z}(-z)^{1-\alpha-\beta}{}_2F_0(1-\alpha,1-\beta;-,z)$, provided that $a=q^\alpha$, $b=q^\beta$ and $x=z/(1-q)$.
On the other hand, convergent series ${}_1F_1(\alpha;\alpha+1-\beta;1/z)$ and
${}_1F_1(\beta;\beta+1-\alpha;1/z)$ appear in the limit $q\to1^-$ of the right-hand side of~\eqref{moth1}.
The aim of this section is to prove the following theorem.
\begin{theorem}
\label{limmo}
The limit $q\to1^-$ of the new connection formula formally gives the following asymptotic formula
\begin{gather*}
e^{1/z}(-z)^{1-\alpha-\beta}{}_2F_0(1-\alpha,1-\beta;-,z)=\frac{\Gamma(\beta-\alpha)}{\Gamma(1-\alpha)}
(-z)^{-\alpha}{}_1F_1\left(\alpha;\alpha+1-\beta;\frac{1}{z}\right)
\\
\qquad
+\frac{\Gamma(\alpha-\beta)}{\Gamma(1-\beta)}(-z)^{-\beta}{}_1F_1\left(\beta;\beta+1-\alpha;\frac{1}{z}
\right).
\end{gather*}
\end{theorem}

In~\cite{Z0}, Zhang has shown a limit of theta functions, taking the principal value of the logarithm on $\mathbb{C}^*\setminus (-\infty ,0]$.

\begin{proposition}
For any $\gamma\in\mathbb{C}^{*}$, we have
\begin{gather*}
\lim_{q\to1^-}\frac{\theta\left(q^\gamma\frac{u}{1-q}\right)}{\theta\left(\frac{u}{1-q}\right)}
(1-q)^{-\gamma}=u^{-\gamma}
\end{gather*}
converges uniformly on compact subset of $\mathbb{C}\setminus(-\infty,0]$.
\end{proposition}
We also remind the formulas for the $q$-gamma function $\Gamma_q(\cdot)$ and the $q$-exponential func\-tion~$E_q(\cdot)$.
The $q$-gamma function is def\/ined by
\begin{gather*}
\Gamma_q(x):=\frac{(q;q)_\infty}{(q^x;q)_\infty}(1-q)^{1-x},
\qquad
0<q<1.
\end{gather*}
This function satisf\/ies $\lim\limits_{q\to1^-}\Gamma_q(x)=\Gamma(x)$~\cite{GR}.
The $q$-exponential function
\begin{gather*}
E_q(z)=\sum_{n\ge0}\frac{q^{n(n-1)/2}}{(q;q)_n}z^n=(-z;q)_\infty
\end{gather*}
satisf\/ies the limit
\begin{gather*}
\lim_{q\to1^-}E_q\left(z(1-q)\right)=e^z.
\end{gather*}

We set $a=q^\alpha$, $b=q^\beta$ and $x=z/(1-q)$ in Theorem~\ref{thm1}.
We introduce the constant
\begin{gather*}
w(\alpha,\beta;q):=(q;q)_\infty(1-q)^{1-\alpha-\beta}.
\end{gather*}
Consider the limit when $q\to1^-$ of each side of the identity of Theorem~\ref{thm1}.
The limit of the left hand side of~\eqref{eq16} is given by the following lemma.
\begin{lemma}%\label{lim1}
For any $\alpha,\beta\in\mathbb{C}^{*}$, $\alpha-\beta\not\in\mathbb{Z}$, we have
\begin{gather*}
\lim_{q\to1^-}w(\alpha,\beta;q)\frac{\left(\frac{q^{\alpha+\beta}z}{1-q};q\right)_\infty}
{\theta\left(-\frac{qz}{1-q}\right)}{}_2\varphi_1\left(q^{1-\alpha},q^{1-\beta};0;q,\frac{q^{\alpha+\beta}z}
{1-q}\right)
\\
\qquad
=(-z)^{1-\alpha-\beta}e^{\frac{1}{z}}{}_2F_0(1-\alpha,1-\beta;-,z).
\end{gather*}
\end{lemma}
\begin{proof}
Exploiting the fact
\begin{gather*}
w(\alpha,\beta;q)\frac{\left(\frac{q^{\alpha+\beta}z}{1-q};q\right)_\infty}{\theta\left(-\frac{qz}{1-q}
\right)}{}_2\varphi_1\left(q^{1-\alpha},q^{1-\beta};0;q,\frac{q^{\alpha+\beta}z}{1-q}\right)
\\
\qquad{}
=\left\{\frac{\theta\left(q^{\alpha+\beta}\left(\frac{-z}{1-q}\right)\right)}{\theta\left(\frac{-z}{1-q}
\right)}(1-q)^{-\alpha-\beta}\right\}\left\{\frac{\theta\left(\frac{-z}{1-q}\right)}
{\theta\left(q\left(\frac{-z}{1-q}\right)\right)}(1-q)\right\}
\\
\qquad\quad{}
\times\frac{1}{E_q\left(-\frac{(1-q)}{q^{\alpha+\beta-1}z}\right)}{}_2\varphi_1\left(q^{1-\alpha}
,q^{1-\beta};0;q,\frac{q^{\alpha+\beta}z}{1-q}\right),
\end{gather*}
we obtain the conclusion.
\end{proof}

Consider the right-hand side of~\eqref{eq16}.
\begin{lemma}%\label{lim2}
For any $\alpha,\beta\in\mathbb{C}^{*}$, $(\alpha-\beta\not\in\mathbb{Z})$, we have
\begin{gather*}
\lim_{q\to1^-}w(\alpha,\beta,q){}_2f_1\left(q^\alpha,q^\beta;q,\frac{z}{(1-q)}\right)
\\
\qquad
=\frac{\Gamma(\beta-\alpha)}{\Gamma(1-\alpha)}(-z)^{-\alpha}{}_1F_1\left(\alpha;\alpha+1-\beta;\frac{1}{z}\right)
+\frac{\Gamma(\alpha-\beta)}{\Gamma(1-\beta)}(-z)^{-\beta}{}_1F_1\left(\beta;\beta+1-\alpha;\frac{1}{z}
\right).
\end{gather*}
\end{lemma}
\begin{proof}
Noting that
\begin{gather*}
w(\alpha,\beta;q)\frac{\left(q^{1-\alpha};q\right)_\infty}{(q^{\beta-\alpha},q;q)_\infty}
\frac{\theta\left(-\frac{q^{\alpha+1}z}{1-q}\right)}{\theta\left(-\frac{qz}{1-q}\right)}\, {}
_2\varphi_1\left(q^\alpha,0;q^{\alpha+1-\beta};q,\frac{q^{1-\alpha-\beta}(1-q)}{z}\right)
\\
\qquad{}
=\left\{\frac{(q^{1-\alpha};q)_\infty}{(q;q)_\infty}(1-q)^\alpha\right\}\left\{\frac{(q;q)_\infty}
{(q^{\beta-\alpha};q)_\infty}(1-q)^{-(1-\beta+\alpha)}\right\}
\\
\qquad\quad{}
\times\left\{\frac{\theta\left(q^{\alpha+1}\left(\frac{-z}{1-q}\right)\right)}{\theta\left(\frac{-z}{1-q}
\right)}(1-q)^{-\alpha-1}\right\}\left\{\frac{\theta\left(\frac{-z}{1-q}\right)}
{\theta\left(q\left(\frac{-z}{1-q}\right)\right)}(1-q)\right\}
\\
\qquad\quad{}
\times{}_2\varphi_1\left(q^\alpha,0;q^{\alpha+1-\beta};q,\frac{q^{1-\alpha-\beta}(1-q)}{z}\right),
\end{gather*}
we prove the lemma.
\end{proof}
Finally, we obtain the proof of Theorem~\ref{limmo}.

\subsection*{Acknowledgements}
The author would like to give heartfelt thanks to Professor Yousuke Ohyama
who provided carefully considered feedback and many valuable comments.
The author also would like to thank the anonymous referees for their helpful comments.

\pdfbookmark[1]{References}{ref}
\LastPageEnding

\end{document}